\documentclass[a4paper,12pt,oneside]{article}
\usepackage[top=3cm, bottom=3cm, left=3cm, right=3cm,dvipdfm]{geometry}
\usepackage{amsmath,amssymb,graphicx}
\usepackage[dvipdfm]{hyperref}
\usepackage{xcolor}
\hypersetup{
    colorlinks,
    linkcolor={blue!50!black},
    citecolor={red!40!black},
    urlcolor={green!60!black}
}

\newtheorem{thm}{\bfseries Theorem}
\newtheorem{lem}[thm]{\bfseries Lemma}        
\newtheorem{cor}[thm]{\bfseries Corollary}     

\newtheorem{cl}[thm]{\bfseries Claim}

\def\rho{\varrho}
\def\qed{\nobreak\hfill $\Box$\par\smallskip}
\def\qq{\nobreak\hfill $\Box$ $\Box$\par\smallskip}

\newcommand{\prf}{\noindent {\bf Proof. }}
\newcommand\prft[1]{\noindent {\bf Proof #1.}}

\newcommand*\compl[1]{V\!-\!({#1})}%
\newcommand*\compli[1]{V\!-({#1})}%
\newcommand*\compls[1]{V\!-\!{#1}}%
\newcommand*\closure[1]{\overline{#1}}%
\newcommand*\Arb{\mathrm{Arb}}%

\title{Spanning tree with lower bound on the degrees} 

\author{Zolt\'an Kir\'aly\\[-0.8ex]
\small E\"otv\"os Lor\'and University \\[-0.8ex]
\small Department of Computer Science  and\\[-0.8ex]
\small MTA-ELTE Egerv\'ary Research Group\\[-0.8ex]
\small P\'azm\'any P\'eter s\'et\'any 1/C \\[-0.8ex]
\small Budapest, Hungary, H-1117\\[-0.8ex]
\small \texttt{kiraly@cs.elte.hu}\\
}

\begin{document}

\thispagestyle{plain}
\maketitle

\begin{abstract}
We concentrate on some recent results of Egawa and Ozeki \cite{EO,EO2}, and He
et al.\ \cite{He}. 
We give shorter proofs and polynomial time
algorithms as well.

We present two new proofs for the sufficient condition for having a spanning
tree with prescribed lower bounds on the degrees, achieved recently 
by Egawa and Ozeki
\cite{EO}. 
The first one is a natural proof using induction, and the second
one is a simple reduction to the theorem of Lov\'asz \cite{L}. 
Using an algorithm
of Frank \cite{F} we show that the condition 
of the theorem can be checked in time $O(m\sqrt{n})$, and moreover, in the same
running time -- if the condition is satisfied -- we can also generate the
spanning tree required.  This gives the first polynomial time algorithm for
this problem.  

Next we show a nice application of this theorem for the simplest case of the
Weak Nine Dragon Tree Conjecture, and for
the game
coloring number of planar graphs, 
first discovered by He et al.\ \cite{He}.

Finally, we give a shorter proof and a polynomial time
algorithm for a good characterization of having a spanning tree with
prescribed degree lower bounds,  for the special case when $G[S]$ is a cograph,
where $S$ is the set of the vertices having degree lower bound prescription  
at least two. This theorem was 
proved by Egawa and Ozeki \cite{EO2} in 2014 while 
they did not give a polynomial
time algorithm. 
 \end{abstract}

\section{Introduction}

Let $G=(V,E)$ be a simple undirected graph, $S\subseteq V$ and $f: S\to
\{2,3,4,\ldots\}$ be an integer-valued function on $S$.
For a subset $X$ of vertices let $f(X)=\sum_{x\in X} f(x)$. 
For disjoint sets of vertices $X$ and $Y$, $d_G(X,Y)$ denotes the number of
edges between $X$ and $Y$, $d_G(X)=d_G(X, \compls{X})$ and $d_G(u)=d_G(\{u\})$.
When the graph $G$ is clear from the context, we omit it from the notation.

The open
neighborhood is denoted by $\Gamma_G(X)=\{u\in \compls{X} \;|\; \exists x\in X,
\; ux\in E\}$, and the closed neighborhood is denoted by
$\Gamma^*_G(X)=\Gamma_G(X)\cup X$. A subgraph induced by a vertex set
$X\subseteq V$ is denoted by $G[X]$,  the number of its edges by $i_G(X)$,
and the number of its components by $c(G[X])$ or  $c_G(X)$. 
We will use the convention that
$\Gamma_G(\emptyset)=\emptyset$ and $c(G[\emptyset])=0$.

Egawa and Ozeki proved the following sufficient condition for having a forest
(or spanning tree) with prescribed lower bounds on the degrees.

\begin{thm}[\cite{EO}]\label{th1}
  If for all nonempty subsets $X\subseteq S$ we have $|\Gamma^*_G(X)|>
  f(X)$ then there is a forest subgraph $F$ of \ $G$, such that for all
  vertices $v\in S$ we have $d_F(v)\ge f(v)$.
\end{thm}

\begin{cor}[\cite{EO}]\label{cor1}
  If for all nonempty subsets $X\subseteq S$ we have $|\Gamma^*_G(X)|>
  f(X)$ and $G$   is connected, then
there is a spanning tree $T$ of $G$, such that for all vertices $v\in S$ we
have $d_T(v)\ge f(v)$.
\end{cor}

Special cases of this theorem appeared in the literature as follows.
When $G$ is bipartite and $S$ is one of the classes, it was proved by Lov\'asz
in 1970 \cite{L}. For general $G$, if $S$ is a stable set, 
it was proved by Frank in 1976 \cite{F} in a stronger form, as in this case
the condition above is also necessary, giving a special case of Theorem
\ref{thgc}
in Section \ref{sec:charact}. For a not necessary stable $S$ a
stronger condition is proved to be sufficient by Singh and Lau \cite{SL}, 
namely: $|\Gamma^*_G(X)|> f(X) + c_G(X)$.

Deciding whether for a triplet $(G,S,f)$ there is a spanning tree $T$
with degree lower bounds, i.e., $d_T(v)\ge f(v)$ for all $v\in S$, is
NP-complete (let $S=V-\{u,v\}$ and  $f(x)=2$ for each $x\in S$; any
appropriate spanning tree is a Hamiltonian path).
However, consider the following algorithmic problem. For given $(G,S,f)$ check,
whether the condition of Corollary \ref{cor1} is satisfied, and if yes, then
construct the appropriate spanning tree $T$. We show that this problem is
polynomially solvable, namely in time $O(m\sqrt{n})$, where $n=|V|$ and $m=|E|$.

In the next section we give a simpler proof than that of Egawa and Ozeki, using
induction. In Section \ref{sec:pr2} we give another proof, that is a simple
reduction to the theorem of Lov\'asz, yielding also a fast algorithm, detailed
in Section \ref{sec:alg}.  In Section \ref{sec:appl} we show an
application (as an example) of Theorem \ref{th1} for the game coloring
number of some planar graphs.
Finally, in Section \ref{sec:charact} we show how we can use these ideas to
prove a good characterization of Egawa and Ozeki \cite{EO2} for a special
case. Our proof is not only shorter but also yields the first polynomial
time algorithm for this case.

\section{First proof -- by induction}\label{sec:pr1}

We prove Theorem \ref{th1} by induction on the number of edges.
If $G$ is a forest or $S=\emptyset$ then the theorem is obviously true.

We call a set $X\subseteq S$ tight if it satisfies the condition
$|\Gamma^*_G(X)|\ge
  f(X)+1$
with equality.

If $uv$ is an edge and $G-uv$ satisfies the condition then we are done by
induction.
So we may assume that for every edge $uv$ the graph  $G-uv$ has a set 
$X\subseteq S$ violating the condition (a violating set). 
This implies that there are no edges outside $S$, and also that for each edge
$uv$ either $u$ or $v$ is contained in a tight set $X$, where the other one
is connected to $X$ by exactly one edge. 
If $u$ is contained in tight set $X$ with $d(v,X)=1$ then we orient edge
$uv$ from $v$ to $u$, otherwise, from $u$ to $v$.
(If both $u$ and $v$ is contained in such a tight set, we choose arbitrarily.)
This oriented graph $\vec G$ has the property that no arc leaves $S$.
(The word \emph{arc} will always refer to a directed edge, in this section a
directed edge of  $\vec G$.)
The in-degree of a vertex $u$ (set $X$) is denoted by $\varrho(u)$
(or $\varrho(X)$ resp.). 

\begin{cl}\label{cl1}
  For each $u\in S$ we have  $f(u)\ge \varrho(u)$.
\end{cl}

\prft{of the Claim} If $\varrho(u)>0$ then $u$ is contained in a tight set. 
As $|\Gamma^*_G(X)|$ is a submodular set function, the intersection and the
union of two intersecting tight sets are both tight.
Thus the intersection $I(u)$ of all
tight sets containing $u$ is also a tight set.
Every arc   $vu$ of $\vec G$ was oriented this way because it entered a tight
set containing $u$, consequently, it must enter $I(u)$ as well.

If $|I(u)|=1$ then, by the tightness, we have $f(u)=d_G(u)\ge \varrho(u)$.
Otherwise, as $I(u)-u$ is not a violating set, if $vu$ is an arc of 
$\vec G$, then  the vertex $v$
does  not have any neighbors in $I(u)-u$. Thus  we have 
$f(I(u))+1-f(u)=f(I(u)-u)+1\le 
|\Gamma^*_G(I(u)-u)|\le
 |\Gamma^*_G(I(u))|-\varrho(u)=f(I(u))+1-\varrho(u)$, giving the claim.\qed

To finish the proof of the theorem it is enough to prove that $G$ is a forest. 
Suppose this is not the case. Choose a cycle $C$ which minimizes $|V(C)-S|$. 
Let $X=V(C)\cap S$ and let $\closure X$ be the closure of $X$ relative to $S$:
$\closure X=\{v\in S \;|\; \exists
x\in X,$ such that $v$ and $x$ are in the same component of $G[S]\}$. 
Clearly $ c(G[\closure X])\le c_G(X)$ and, by the observation made above, no arc
leaves $\closure X$.

If $V(C)\subseteq S$ then, using  Claim \ref{cl1} and the fact
that $G[\closure X]$ is now connected and contains a cycle,
$f(\closure X)\ge i_G(\closure X)+\varrho(\closure X)\ge |\closure
X|+\varrho(\closure X) 
\ge |\Gamma^*_G(\closure X)|$, and this contradicts to the assumption of the
theorem. 

Otherwise, $G[S]$ is a forest and $|V(C)-S|\ge  c(G[\closure X])$.
Now, by Claim \ref{cl1},
$f(\closure X)\ge i_G(\closure X)+\varrho(\closure X)\ge  |\closure
X|-c(G[\closure X])+ 
\varrho(\closure X)$.
As $V\!-\!S$ is an independent set and no arc leaves $S$, 
at least two arcs go from any vertex of $V(C)-S$ to $\closure X$,
that is $\closure X$ has at most 
$\varrho(\closure X)-|V(C)-S|\le \varrho(\closure X)-c(G[\closure X])$ 
different neighbors in $V\!-\!S$.
Thus $|\Gamma^*_G(\closure X)|\le |\closure X|+\varrho(\closure X)-c(G[\closure
X])\le  
f(\closure X)$, a contradiction again. \qq

\section{Second proof -- reduction to Lov\'asz' theorem}
\label{sec:pr2}

In this section we prove Theorem \ref{th1} using a theorem of Lov\'asz
\cite{L}. We quote this old theorem reformulated for fitting the notation used
in this paper. We denote by $f^+$
the function $f+1$, i.e., $f^+(x')=f(x')+1$ for $x'\in S'$.

\begin{thm}[\cite{L}]\label{th2}
  Let $B=(S'\cup V,E')$ be a bipartite graph and $f: S'\to
\{2,3,4,\ldots\}$ be a function. $B$ has a forest subgraph $F_0$ with the
property $d_{F_0}(x')=f^+(x')$ for every $x'\in S'$, if and only if for all
nonempty $X'\subseteq S'$ we have $|\Gamma_B(X')|> f(X')$.
\end{thm}

\prft{of Theorem \ref{th1}} We have $(G, S, f)$ given, and let $S'$ be a set
disjoint from $V$ with elements $S'=\{u' \;|\; u\in S\}$, and extend $f$ to
$S'$ in the obvious way: $f(u'):=f(u)$ for each $u\in S$.

Construct a bipartite graph $B=(V\cup S', E')$ as
follows. For each ordered vertex pair $(u\in S,\; v\in V)$ we put an edge $u'v$
into $E'$ if $uv\in E$,  
and we also put the 'vertical' edges $u'u$ for each $u\in S$.
Observe that for each $X\subseteq S$ (if $X'$ denotes the corresponding subset
of $S'$) we have $\Gamma_B(X')=\Gamma^*_G(X)$. Therefore the condition of
Theorem \ref{th2} is satisfied and thus we have a forest subgraph $F_0$
of $B$ with
$d_{F_0}(x')=f^+(x')$ for every $x'\in S'$. 

First we claim that we can modify $F_0$ to get another forest subgraph $F_1$
of $B$, so that $d_{F_1}(x')\ge f^+(x')$ for every $x'\in S'$,
with the additional property that $F_1$ contains every vertical edge.
Suppose $u'u\not\in E'(F_0)$. If $u'$ and $u$ are in different components of
$F_0$, then we add the edge $u'u$. Otherwise, there is unique path
$u'vu_1'\ldots u_t'u$ in $F_0$, in this case we delete the
edge $u'v$ and add the edge
$u'u$, still resulting in a forest (with the same degrees inside $S'$).

\begin{figure}[!ht]
      \centering
      \smallskip
      \includegraphics{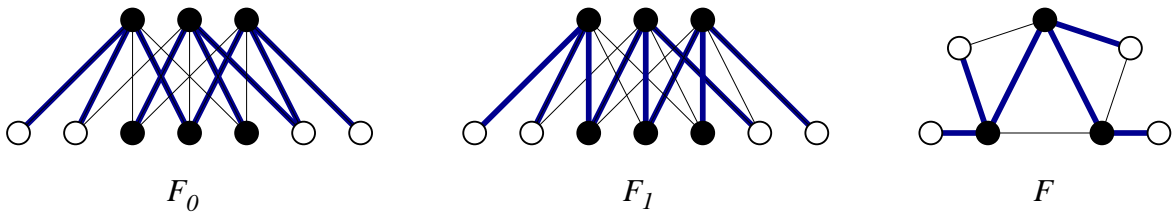}
      \caption{An example for the proof, $F_0, F_1$ and $F$ are shown in blue
        and bold.}\label{fig1}
\end{figure}

Finally, we construct the desired forest $F$ by contracting each vertical edge
(we contract $u'$ to vertex $u$). It is easy to see, that in this way $F$
becomes a forest subgraph of the graph $G$, and $d_F(x)\ge d_{F_1}(x')-1\ge
f(x)$ for  every $x\in S$. 
See Figure \ref{fig1} for an example, where vertices if $S\cup S'$ are black,
and $f(u)=2$ for every $u\in S$. \qed
 
{\sl Remark: Theorem \ref{th2} remains true if we allow $f(x')=1$ for some
  $x'\in S'$.}

\section{A polynomial time algorithm for checking the condition
and  constructing  the tree}\label{sec:alg}

In this section we first describe the algorithmic proof of Frank \cite{F} for
Theorem \ref{th2}. After cloning each vertex $u'\in S'$ into $f(u')$ copies
and running e.g., the algorithm of Hopcroft and Karp 
for maximum bipartite matching,
we either get a forest $F'$ with degrees $d_{F'}(u')=f(u')$ for every $u'\in
S'$ (and $d_{F'}(v) \le 1$ for each $v\in V$), or we get a subset $X'\subseteq
S'$ that violates the Hall condition, namely $|\Gamma_B(X')|< f(X')$. In this
latter case the corresponding $X\subseteq
S$ clearly violates the condition of the theorem as well.

\begin{figure}[!ht]
      \centering
      \smallskip
      \includegraphics{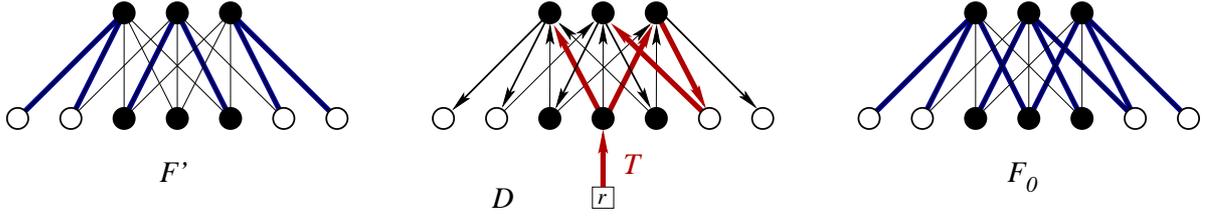}
      \caption{An example for Frank's algorithm, $F'$ and $F$ are shown in blue
        and bold, $T$ is shown in red and bold.}\label{fig2}
\end{figure}

We make an auxiliary digraph $D=(U,A)$, where $U=V\cup S'\cup \{r\}$.
We orient edges of $F'$ from $S'$ to $V$, other edges of $B$ from $V$ to $S'$,
and finally add arcs $rv$ for each $v\in V$ uncovered by $F'$. 
We run a BFS from vertex $r$ in digraph $D$. This gives an arborescence $T$
rooted at $r$ which spans all vertices  reachable from $r$.
If every vertex in $S'$ is reachable from $r$, then for each arc $vu'$ of $T$
leading from $V$ to $S'$ we add the corresponding edge to $F'$ resulting in the
desired forest $F_0$ in $B$ (these are not edges of $F'$, so they increase the
degree of every $u'\in S'$). Observe that we did not create any cycle because 
every vertex $u'\in S'$ has in-degree one in $T$ and the arborescence
$T$ does not contain any directed cycle.)
See Figure \ref{fig2} for an example.   

Otherwise, if $X'$ denotes the set of vertices
of $S'$ that are not reachable from $r$, then we claim that $X'$ violates the
condition of Theorem \ref{th2}. If not, then there exist
 $u\in V$ and $x'\in X'$ such that 
$ux'\in A$ but either $u$ is uncovered by $F'$ or 
$u$ is a leaf of $F'$, and its unique $F'$-neighbor $y'$ is
in $S'\!-\!X'$. 
In both cases $u$ is reachable from $r$
(in the first case $ru$ is an arc, in the second case $y'$
is reachable from $r$ and $y'u$ is an arc);
consequently, $x'$ is also
reachable from $r$, a contradiction.

\paragraph{Final algorithm and running time} 
We are ready to give an 
algorithm running in time $O(m\sqrt{n})$
for deciding whether the condition of 
Corollary \ref{cor1} is satisfied or not. Moreover,  
-- if the condition holds -- we can also generate the spanning tree required
in the same running time.  This gives the first polynomial time algorithm for
this problem.

First we construct the bipartite graph $B$ in time $O(m)$.
Then we follow the steps of Frank's algorithm.
Observe that for running the algorithm of Hopcroft and Karp, we
do not need to make the cloning in reality; it is enough to do it imaginarily.
Doing so keeps the running time $O(m\sqrt{n})$. 
Constructing $D$, running BFS and constructing $F_0$ can be done in
time $O(m)$.

Next we make $F_1$ from $F_0$ and then we construct $F$ by contracting the
vertical edges
as in the proof presented in Section \ref{sec:pr2}, 
these are algorithmically easy jobs,
they can be done in time $O(m)$.

Finally we make the desired spanning tree from forest $F$, 
it  can also be done in time $O(m)$.

\section{Some applications: WNDT Conjecture and \\
game coloring number of 
planar graphs}\label{sec:appl}

We show an interesting application of Theorem \ref{th1} as an example.  
For a subset $X$ of vertices, if $|X|>1$, then we define $\lambda_G(X)=
\frac{i(X)}{|X|-1}$ and $\Arb(G)=\max \{\lambda_G(X) \;|\;
X\subseteq V,\; |X|>1\}$.
The Weak Nine Dragon Tree (WNDT for short) Conjecture is the following:
if for integers $k$ and $d$ we have  $\Arb(G)\le k+\frac{d}{d+k+1}$, then
there are $k$ forests $F_1,\ldots,F_k$, such that the maximum degree in
$G-F_1-\ldots-F_k$ is at most $d$.
The conjecture was proved by Kim et al.\ \cite{Kim} for the case of $d>k$.

Here we show that if we further restrict ourselves to the special case of
$k=1$, then this results in a simple consequence of Theorem \ref{th1}.

\begin{thm}\cite{Kim}\label{th_WNDT_k1}
If $d\ge 2$ is an integer and $\Arb(G)\le 1+\frac{d}{d+2}$, then there is a
forest subgraph $F$ of $G$, such that for every vertex $v$ we have
$d_G(v)-d_F(v)\le d$. 
\end{thm}

Actually we prove a stronger form (also proved in \cite{Kim}).

\begin{thm}\cite{Kim}\label{th_WNDT_k2}
If $d\ge 2$ is an integer and for each nonempty subset $X$ of the vertices we
have $2(d+1)\cdot|X|> (d+2)\cdot i(X)$, then there is a forest
subgraph $F$ of $G$, such that for every vertex $v$ we have $d_G(v)-d_F(v)\le
d$. 
\end{thm}

\prf 
Let $S=\{v\in V\; | \; d(v)\ge d+2\}$ and let $f$ be defined on $S$ by
$f(v)=d(v)-d$. 
For a subset $X\subseteq S$ let $\Gamma_j(X)=\{v\in \compls{X}\; | \;
d(v,X)=j\}$, and let  $\tilde X=X\cup
\bigcup_{j=2}^{|X|} \Gamma_j(X)$.
By the condition of the theorem $(d+1)\cdot|\tilde X|> 
\frac{d+2}2\cdot i(\tilde X)$,
i.e.,  $(d+1)\cdot|X| + (d+1)\cdot\sum_{j=2}^{|X|} |\Gamma_j(X)|
>\frac{d+2}2\cdot \bigl(i(X)+\sum_{j=2}^{|X|} j\cdot|\Gamma_j(X)|\bigr)$. 
Realigned we get
\vspace{-12pt}
\[(d+1)\cdot|X|
>\frac{d+2}2\cdot i(X) +\sum_{j=2}^{|X|}
\Bigl[\bigl(j\cdot\frac{d+2}{2}-(d+1)\bigr)\cdot|\Gamma_j(X)|\Bigr]\ge 
\]
\vspace{-12pt}
\[\ge 2\cdot i(X)+\sum_{j=2}^{|X|}\Bigl[
\bigl((j-2)\cdot\frac{d}{2}+j-1\bigr)\cdot|\Gamma_j(X)|\Bigr]\ge 
2\cdot i(X) +\sum_{j=2}^{|X|}\Bigl[
(j-1)\cdot|\Gamma_j(X)|\Bigr],\]
\vspace{-2pt}

 as $d\ge 2$. We have
$f(X)=2\cdot i(X)+ \sum_{j=1}^{|X|} j\cdot|\Gamma_j(X)| - d\cdot |X|=
2\cdot i(X)+ \sum_{j=2}^{|X|} \Bigl[(j-1)\cdot|\Gamma_j(X)|\Bigr] +
\bigl(|\Gamma(X)|-d\cdot|X|\bigr)<
(d+1)\cdot|X| +\bigl(|\Gamma(X)|-d\cdot|X|\bigr)=|\Gamma^*(X)|$,
thus the condition of Theorem \ref{th1} is satisfied and the forest $F$
produced fulfills the statement of  our theorem.
\qed

Of course, we can apply the algorithm described in the previous section and
efficiently make this decomposition.

Let $G$ be a simple connected planar graph with girth $g\ge 5$.
We know by Euler's formula that $i(X)<\frac{g}{g-2}\cdot |X|$ for every
subset $X$ of the vertices, and $\frac{g}{g-2}\le\frac{2d+2}{d+2}$ if
$d\ge\frac{4}{g-4}$.
We get the following corollary (which is a
strengthening of a theorem proved first by He
et al.\ in \cite{He},  the improvement was reported to
be proved in  \cite{Kleitman}).

\begin{cor}\cite{He,Kleitman}\label{cor_planar}
  If $G$ is a simple connected planar graph with girth  at least $g$ (where
$g=5$ or $g=6$), then
there is a spanning tree $T$ of $G$, 
such that for every vertex $v$ we have $d_G(v)-d_T(v)\le
\frac{4}{g-4}$.
\end{cor}

The game coloring number was defined by Zhu \cite{Zhu} via a two-person game
(for upper bounding the so-called ``game chromatic number'').
Alice and Bob remove vertices of $G$ in turns. The back-degree of a vertex 
is the number of its previously removed neighbors. The game coloring number
$\mathrm{col}_g(G)$ is the smallest $k+1$, where Alice can achieve that every
vertex has back-degree at most $k$. An easy observation of Zhu \cite{Zhu} states
that if the edges of $G$ can be partitioned into graphs $G'$ and $H$, 
then  $\mathrm{col}_g(G)\le \mathrm{col}_g(G') + \Delta(H)$, where $\Delta$
denotes the maximum degree. Faigle et al.\ \cite{Faigle} proved that the game
coloring number of a tree is at most $4$. Consequently, we get the following
result,  that is also a strengthening of a theorem proved by He
et al.\ in \cite{He}). We also note, that by our algorithmic results we also
provide a simple polynomial time algorithm for Alice for winning the game,
as the proofs in \cite{Faigle} and 
\cite{Zhu} are algorithmic.

\begin{cor}\cite{He,Kleitman}\label{cor_planar-col}
  If $G$ is a simple planar graph with girth at least $5$, then
$\mathrm{col}_g(G)\le 8$.
  If $G$ is a simple planar graph with girth at least $6$, then
$\mathrm{col}_g(G)\le 6$.
 \end{cor}

\section{Good characterization for a special case}\label{sec:charact}

In \cite{EO2} Egawa and Ozeki proved the following theorem stating a good
characterization if $G[S]$ is a cograph, i.e., it does not contain an induced
$P_4$. 
By the definition, an induced subgraph of a cograph is a cograph, and for any
two different vertices of the same connected component of a cograph, they are
either adjacent or have a common neighbor. The latter property
 is equivalent to saying that
every component has diameter at most 2.

Egawa and Ozeki also showed by a simple example, that this
characterization does not remain true if $G[S]=P_4$: let the vertices of the
$P_4$ be $v_1, v_2, v_3, v_4$ and let $G$ have two more vertices, $a$ and $b$,
such that $a$ is connected to $v_1$ and $v_4$ while $b$ connected to $v_2$ and
$v_3$; and let $f(v_1)=f(v_2)=f(v_3)=f(v_4)=2$. 

\begin{thm}[\cite{EO2}]\label{thgc}
  If $G[S]$ is a cograph, then $G$ has a forest subgraph with
  degree lower bounds $f: S\to\{2,3,4,\ldots\}$  on $S$ 
if and only if for all nonempty subsets
  $X\subseteq S$ we have
\[
|\Gamma_G(X)|+2|X|-c_G(X)> f(X).
\]
\end{thm}

\prf We follow the outline of the proof in \cite{EO2} but we make some
simplifications resulting in a significantly shorter proof. 
We also give a polynomial time algorithm for finding the appropriate forest.

It is not hard to see that the condition above is necessary
(even for the case when $G[S]$ is not a cograph). Let $F$ be a
forest with $d_F(u)\ge f(u)$ for each $u\in S$, and let $X\subseteq S$. 
Now $i_G(F[X])=|X|-c(F[X])$ and $f(X)\le 2i_G(F[X])+d_F(X)= 2|X|-2c(F[X])+
d_F(X)\le (d_F(X)-c(F[X]))+2|X|-c_G(X)$ and $d_F(X)-c(F[X])\le
|\Gamma_F(X)|-1\le |\Gamma_G(X)|-1$, because $F$ is a forest subgraph of $G$.

For $Z,X\subseteq V$, let $\Gamma_{Z}(X)=\Gamma(X)\cap Z$.
(For notational symmetry we will use the notation also for $\Gamma_{Z}(Z)$,
though this set is always empty.) 
We denote the cograph $G[S]$ by $H$. 

\begin{cl}\label{cl_subm1}
If $A,B\subseteq S$, then
\[
|\Gamma_{A\cup B}(A\cup B)|-c_H(A\cup B)+|\Gamma_{A\cup B}(A\cap B)|-c_H(A\cap
B)\le 
\]
\[\le |\Gamma_{A\cup B}(A)|-c_H(A)+|\Gamma_{A\cup B}(B)|-c_H(B).
\]
\end{cl}

\prf We first prove the claim for the case when $G[A\cup B]$ is a connected
cograph and $A\cap B\ne\emptyset$.
 We use the well-known observation of Erd\H os and Rado stating that a
graph or its complement is connected. 
As  for any $x,y\in A\cap
B$ they are either in the same component of $G[A]$ or in the same component of
$G[B]$, (if they are not connected, then they have a common neighbor in $A\cup
B$), we may assume that $G[A\cap B]$ is inside a component $K$ of  $G[A]$.
Let $K_1,\ldots,K_a$ denote the other components of $G[A]$, and
$L_1,\ldots,L_b,I_1,\ldots,I_c$ denote the components of $G[B]$,
where $I_j$ are the components intersecting $A\cap B$.
As $c_H(A\cap B)\ge c$, it is enough to prove
\[
0+|\Gamma_{A\cup B}(A\cap B)|-1-c\le |\Gamma_{A\cup B}(A)|+|\Gamma_{A\cup B}(B)|
-(a+1)-(b+c),
\]
i.e., $|\Gamma_{A\cup B}(A)|+|\Gamma_{A\cup  B}(B)|-
|\Gamma_{A\cup B}(A\cap B)|\ge a+b$.
This can be easily seen, as
 $\Gamma_{A\cup B}(A\cap B)=(A\cup B)-(\bigcup_{i=1}^a K_i \cup
 \bigcup_{j=1}^b L_j)$ and -- as we assumed that $G[A\cup B]$ is connected --
every $K_i$ contains a vertex connected to $B\!-\!A$ and
every $L_j$ contains a vertex connected to $A\!-\!B$.

If $A\cap B=\emptyset$, then the situation is very similar, we need to  prove
$|\Gamma_{A\cup B}(A)|+|\Gamma_{A\cup  B}(B)|\ge a+b-1$, however, as we
assumed $G[A\cup B]$ to be connected, this is always satisfied with a strict
inequality. 

If $G[A\cup B]$ has several components, then it is enough to 
prove the claim for every component separately, thus the same proof works. 
We remark that for the case $A\cap B=\emptyset$ we still have strict
inequality if there is an edge between $A$ and $B$.
\qed

Let $b_0(X)$ and $b(X)$ be set-functions on the subsets of $S$ defined by
$b_0(X)=|\Gamma(X)|-c_H(X)$, and $b(X)=b_0(X)+2|X|-f(X)$.
Let $A,B\subseteq S$, and denote by $U(A,B)$ the set of vertices in
$\compl{A\cup B}$ connected to both $A\!-\!B$ and $B\!-\!A$ but not to $A\cap B$
(in other words $U(A,B)=\Gamma(A)\cap \Gamma(B)-\Gamma(A\cap B)$).
We claim that $b_0$ and $b$ are submodular, moreover,
$b_0(A\cup B)+b_0(A\cap B)+|U(A,B)|\le b_0(A)+b_0(B)$.
As $b_0(A)=|\Gamma_{A\cup B}(A)|+|\Gamma_{\compli{A\cup B}}(A)|-c_H(A)$,
using Claim \ref{cl_subm1}, it is enough to prove that
\[|\Gamma_{\compli{A\cup B}}(A\cup B)|+|\Gamma_{\compli{A\cup B}}(A\cap B)|
+|U(A,B)|\le |\Gamma_{\compli{A\cup B}}(A)|+|\Gamma_{\compli{A\cup B}}(B)|.\]
However, this is obvious by the definition of $U(A,B)$.
As $b(X)$ is the sum of $b_0(X)$ and the modular function $2|X|-f(X)$, the
same statement holds for $b$ as well.

We call a nonempty subset $X\subseteq S$ {\bf tight} if $b(X)=1$. 

\begin{cor}\label{cor_tight}
Suppose the condition of the
theorem, i.e.,  $b(X)\ge 1$ for all $\emptyset\ne X\subseteq S$ holds.
The intersection and union 
of two intersecting tight sets $A$ and $B$ is tight, and
$|U(A,B)|=0$. 
If $A$ and $B$ are disjoint tight sets, and 
either $U(A,B)\ne\emptyset$ or there is an
edge connecting $A$ and $B$, then $A\cup B$ is tight, moreover either 
$|U(A,B)|=1$ and no edge connects $A$ to $B$; or $U(A,B)=\emptyset$.
\end{cor}

Let $W$ denote the set $\compls{S}$.
We may assume that every vertex $v\in S$ is contained in a tight set,
otherwise, we can increase $f(v)$ without violating the condition.
By Corollary \ref{cor_tight}, $I(v)$, the \emph{intersection of all tight sets
containing} $v$ is also a tight set.
We also suppose that for every edge $wv$ (where $w\in W$) the graph
$G-wv$ would violate the condition, i.e., $v\in S$ and $d(w,I(v))=1$.
(If this is not the case, then we simply
delete the edge $wv$.)
Unfortunately, we may not assume similar condition
about edges induced by $S$ because deleting such an edge can introduce an
induced $P_4$.

We denote the components of $G[S]$ by $Z_1,\ldots,Z_t$.
Call a vertex $u\in Z_i$ {\bf proper} if $I(u)\subseteq Z_i$. Our 
first goal is to
prove that every vertex in $S$ is proper. We need some preliminary
observations.

\begin{cl}\label{cl_aux}
Suppose tight sets $A$ and $B$ intersect the same component $Z_i$ of $G[S]$. 
Then $A\cup B$ is tight.
Consequently, if $u, v\in Z_i$ and $u\ne v$,
then $I(u)\cup I(v)$ is tight, 
moreover,  a vertex $w\in W$ cannot be connected to both $u$ and $v$.
\end{cl}

\prf 
Either the sets $A$ and $B$ are intersecting, or connected by an edge, or
otherwise -- using that $G[Z_i]$ is a cograph -- they have a common neighbor
$x\in Z_i$ which is not in $A\cup B$, so $x\in U(A,B)$.  Thus the
first statement is a consequence of Corollary \ref{cor_tight}.

Suppose $wu$ and $wv$ are edges.
As $wu$ does not enter $I(v)$ (because $wv$ is the unique
edge entering $I(v)$), and $wv$ does not enter $I(u)$, we have 
$w\in U(I(u), I(v))$, so using Corollary \ref{cor_tight} again, we get a
contradiction. 
\qed

\begin{lem}\label{lem_proper}
  Every vertex $u\in S$ is proper.
\end{lem}

\prf 
Suppose this is not the case and let $u\in S$ be an unproper vertex for which
$I(u)$ is minimal. Let $I=I(u)$, and denote the components of $G[S]$ 
intersected by $I$ by  $Z_1,\ldots,Z_r$.
Take any vertex $v\in I$. If $I(v)\ne I$, then  $I(v)\subset I$ 
by Corollary \ref{cor_tight}, and thus $v$ is proper by the minimality of $I$.

Let $E''=\{ wv\in E\;|\; w\in \Gamma_W(I),\; v\in I\}$. 
For $i=1\ldots r$, let $A_i=I\cap Z_i$.
   
\noindent
First we claim that if we consider the subgraph defined by $E''$, 
and thereafter we
contract each $A_i$ to a vertex $a_i$, 
then the resulting bipartite graph is a
forest. 

Suppose not, i.e., it
contains a cycle, wlog 
$C'=w_1,a_1,w_2,a_2,\ldots,w_k,a_k,w_1$ in cyclic order.
Let $C$ be its ``pre-image'' in $G$, a cycle 
\[w_1,v_1,x_1,u_2,w_2,v_2,x_2,u_3,w_3,\ldots,u_k,w_k,v_k,x_k,u_{k+1}=u_1,w_1,\]
where $w_i\in \Gamma_W(I)$,  $v_i,u_{i+1}\in A_i$ and $x_i\in Z_i$ (the vertices
$v_i$ and $u_{i+1}$ are well defined by the edges of $C'$, 
and we connect them with a shortest
path inside $Z_i$).
Note that $v_i,x_i,u_{i+1}$ are
not necessarily distinct vertices, so two subsequent vertices of this sequence
are either identical or connected by an edge. 
As every $w_i$ is connected to two contracted vertices, both $v_i$ and
$u_{i+1}$ are proper vertices, otherwise, e.g., 
$I(v_i)=I$ and the edge $w_iv_i$ is not a
unique edge from $w_i$ that enters $I(v_i)$. 

Using  Claim \ref{cl_aux}, the sets
$B_i=I(v_i)\cup I(u_{i+1})\subseteq A_i$ are tight sets.  
By repeatedly using Corollary \ref{cor_tight} we get that $D_j=\cup_{i=1}^j
B_i$ are tight sets for $j=1,2,\ldots,k-1$; note that $D_{k-1}\subseteq
A_1\cup\ldots\cup A_{k-1}$. Finally, we get a contradiction to Corollary 
\ref{cor_tight}
for disjoint tight sets $D_{k-1}$ and $B_k$ as $w_k,w_1\in U(D_{k-1},B_k)$. 

Now we are able to finish the proof of Lemma \ref{lem_proper}.
By our assumptions, $b(A_i)\ge 1$ for $i=1\ldots r$, and there is a $j$ such
that $u\in A_j$, we have $b(A_j)\ge 2$ as $A_j$ is not tight. 
Thus we have $\sum_{i=1}^r b(A_i)\ge r+1$.  
We claim that $|\Gamma_W(I)|+r\le |\Gamma_W(I)|+\sum_{i=1}^r b(A_i) -
b(I)=|E''|$,  
this would give
the required contradiction.
We have $b(I)=1$ because $I$ is tight. 
The sets $\Gamma_S(A_i)\subseteq Z_i$ are pairwise
disjoint and $I=\cup_{i=1}^{r} A_i$; therefore $\Gamma_S(I)$
is the disjoint union of sets $\Gamma_S(A_i)$. 
Moreover $c_G(I)=\sum_{i=1}^{r} c_G(A_i)$.
Hence $\sum_{i=1}^r b(A_i) - b(I)=\sum_{i=1}^{r} |\Gamma_S(A_i)|+
\sum_{i=1}^{r} |\Gamma_W(A_i)|-|\Gamma_S(I)|-|\Gamma_W(I)|=
\sum_{i=1}^{r} |\Gamma_W(A_i)|-|\Gamma_W(I)|=
\sum_{i=1}^{r} d(A_i,\Gamma_W(I))-|\Gamma_W(I)|$ by the second
statement of Claim
\ref{cl_aux},
and $\sum_{i=1}^{r} d(A_i,\Gamma_W(I))=|E''|$, so our last claim is proved,
finishing 
the proof of the lemma.
\qed

For any component $Z_i$ of $G[S]$, we have $Z_i= \cup_{u\in Z_i} I(u)$ 
by Lemma \ref{lem_proper}; thus repeated usage of
Claim \ref{cl_aux} shows that $\cup_{u\in Z_i} I(u)$ is a tight set.

\begin{cor}\label{cor_proper}
  For every component $Z_i$ of $G[S]$, the set $Z_i$ is tight.
\end{cor}

We make an auxiliary bipartite graph $G'$ by contracting each component of
$G[S]$ (we delete the loops arising).
By Claim \ref{cl_aux} no parallel edges can arise.
The components of $G[S]$ are $Z_1,Z_2,\ldots,Z_t$, hence the
corresponding
contracted vertices of $G'$ will be denoted by $z_1,z_2\ldots,z_t$.

We prove Theorem \ref{thgc} by induction on the number of vertices. 
If $G'$ has an isolated vertex,
then this vertex is either $w\in W$ (we simply delete $w$ and use induction),
or a vertex $z_j$. In this latter case $G[Z_j]$ is a component of $G$,
so we can use induction separately 
for $G[Z_j]$ and for $G-Z_j$.

If $G'$ has a vertex of degree one, then it is either $w\in W$ or a vertex
$z_j$. If $w$ has degree one we take its neighbor $u\in S$, delete $w$ and
reset $f(u)=f(u)-1$. Now we can use induction, the assumption of the
theorem is not violated (it may be the case that $u$ gets outside of $S$ 
-- if $f(u)$ becomes $1$  -- but
$G[S-u]$ remains a cograph). Suppose $z_j$ has degree one in $G'$ and $uw$ is
the unique edge leaving $Z_j$ in $G$, where $u\in Z_j$ and $w\in W$. We delete
edge $uw$, reset $f(u)=f(u)-1$, and then we can use induction separately for
$G[Z_j]$ and for $G-Z_j$, finally the edge $uw$ can be put back safely to the
union of the two resulting forests.

Otherwise, we have a cycle in $G'$ with vertices $w_1,z_1,w_2,z_2\ldots,
w_k,z_k$.
We repeat the arguments as in the proof of Lemma \ref{lem_proper} in order
to show that this assumption leads to a contradiction.
Now let $D_j=\cup_{i=1}^j Z_i$, 
these are tight sets for $j=1,2,\ldots,k-1$ by Corollaries
\ref{cor_proper} and \ref{cor_tight},
and finally we get a contradiction
for tight sets $D_{k-1}$ and $Z_k$ as $w_k,w_1\in U(D_{k-1},Z_k)$.
\qq

\subsection{Algorithmic aspects}

Egawa and Ozeki already observed that their proof is ``almost''  algorithmic
but they were not able to give a polynomial time algorithm.
They wrote: `we believe that there is a polynomial time algorithm to find an
$(X,f)$-tree in a graph satisfying condition (1).'

The author of the present paper can only guess at the reason of this.
We think of two possibilities.
Actually, they also used induction 
but they used the induction hypothesis for every
tight set $I\ne S$, for the same graph with
$S'=I$, thus  exponentially many times. However,
they did not really need the forests arising from the hypothesis, only their
existence. 
On the other hand, they wrote: `To find an appropriate vertex or edge ...'
Probably they did not realize that they do not need to find an appropriate
edge. It is enough to check the condition for graph $G-wv$ for every edge $wv$
where $w\in W$ and $v\in S$. If the condition still holds for $G-wv$, then
the edge
$wv$ can be deleted. If for none of the edges it holds, 
then their Claim 8 applies, and so they could make
the recursion for $G-wv$ for \emph{any} edge between $S$ and $W$. 
However, this  train of thought is not obvious at all, one must check their 
long proof thoroughly.

Our proof uses the inductive hypothesis only once, so we may call the
forest-construction procedure for two graphs with total number $|V|$ of the
vertices. 
So from this proof, it is easy to conclude 
that by using general strongly polynomial 
submodular function minimization (SPSFM hereafter) of either Iwata,
Fleischer and Fujishige \cite{Iwata} or of Schrijver \cite{Sch} we have 
a polynomial time algorithm for constructing the desired forest if the
condition of the theorem holds (Egawa and Ozeki already showed that checking
the condition can be done by one call of SPSFM).

To be more precise, our algorithm is as follows.
First we delete edges inside $W$ and check the condition for $G, S, f$.
Then for each $v\in S$ we check whether the condition still holds if we
increase $f(v)$ by one. If yes, then we increase $f(v)$ and continue.
While the condition holds, we have $f(S)<2|V|$, 
so this process can be done by at most $2|V|$ calls of
SPSFM.

When none of the $f$ values can be increased, then
every vertex $v\in S$ contained in a tight set. We claim that we can also get
the sets $I(v)$ for all $v\in S$.
This can be done in many ways, for example minimizing the submodular function
$b'(X)$, where $b'(X)=|V|\cdot b(X)+|X|$ if $v\in X$ and
$b'(X)=3|V|+|V|\cdot b(X)+|X|$ otherwise. In short, after $O(|V|)$ calls of
SPSFM we ensured that every vertex $v\in S$ is in a tight set and we calculated
$I(v)$.

Next we check that for every edge $wv$ (where $w\in W$ and $v\in S$)
$d(w,I(v))=1$ or not. We do not need any further call of SPSFM, this can be
done in $O(|V||E|)$ steps. If $d(w,I(v))>1$, then we delete the edge $wv$.

Finding an isolated vertex or a leaf of $G'$ is easy.
The main point for the remaining part is that we do not need to recalculate
anything when making recursive calls for graphs $G[Z_i]$ and $G-Z_i$ or
when we delete a vertex $w\in W$. This is because the same sets $I(v)$ do the
job. 

In conclusion, we can find the appropriate forest by $O(|V|)$ calls of
 SPSFM and by $O(|V||E|)$ simple graph operations.

\section{Acknowledgment} 

Special thanks to Andr\'as Frank, who, (after I wrote down my first simple
proof and gave the first polynomial time algorithm), 
suggested the bipartite graph
construction used here, in order for simpler checking 
the condition; and who taught me his nice algorithm.

Research is supported by a grant (no.\ K 109240) from the National Development
Agency of Hungary, based on a source from the Research and Technology
Innovation Fund.

\end{document}